\documentclass[10pt,reqno,oneside]{amsproc}

\usepackage{amsfonts}
\usepackage{dsfont}
\usepackage[usenames,dvipsnames,svgnames,table]{xcolor}
\usepackage{amsmath, amsthm, amssymb, enumerate, mathscinet}
\usepackage[T1]{fontenc}

\usepackage{color}  
\usepackage{marginnote}
\usepackage[colorlinks=true, pdfstartview=FitV, linkcolor=black, citecolor=black, urlcolor=black]{hyperref}

  \chardef\forshowkeys=0
  \chardef\refcheck=0
  \chardef\showllabel=0
  \chardef\sketches=0
  \chardef\showcolors=0

\ifnum\forshowkeys=1
  
  \usepackage[notref,notcite,color]{showkeys}
\fi

\ifnum\showllabel=1
\def\llabel#1{\marginnote{\color{gray}\rm(#1)}[-0.0cm]\notag}
\else
\def\llabel#1{\notag}
\fi

 \ifnum\refcheck=1
  \usepackage{refcheck}
\fi

\newtheorem{theorem}{Theorem}[section]

\newtheorem{Lemma}[theorem]{Lemma}

\theoremstyle{definition}
\newtheorem{definition}[theorem]{Definition}
\newtheorem{Remark}[theorem]{Remark}

\def\lec{\lesssim}

   \def\startnewsection#1#2{\section{#1}\label{#2}\setcounter{equation}{0}}

   \def\comma{ {\rm ,\qquad{}} }            
   
\def\div{\mathop{\rm div}\nolimits}
\def\curl{\mathop{\rm curl}\nolimits}

\ifnum\showcolors=1
  \def\cole{\color{coloroftheorems}}
  \definecolor{colorcccc}{rgb}{0.7,0.7,0.7}

  \def\colb{\color{black}}
  \definecolor{colorpppp}{rgb}{0.6,0.0,0.1}
  \definecolor{colorgggg}{rgb}{.0,0.4,0.0}
  \definecolor{colorhhhh}{rgb}{0,0.6,0.2}
  \definecolor{colorgray}{rgb}{0.8,0.8,0.8}

  \definecolor{coloroftheorems}{rgb}{0.6,0.0,0.6}
  \definecolor{colorigor}{rgb}{1, 0.2, 0.8}
  \definecolor{amethyst}{rgb}{0.6, 0.4, 0.8}

  \definecolor{colororange}{rgb}{0.8,0.2,0}
  \definecolor{colorpurple}{rgb}{0.6,0.0,0.6}
  
\else   
  \def\cole{}
  \definecolor{colorcccc}{rgb}{0,0,0}

  \def\colb{\color{black}}
  \definecolor{colorpppp}{rgb}{0,0,0}
  \definecolor{colorgggg}{rgb}{0,0,0}
  \definecolor{colorhhhh}{rgb}{0,0,0}
  \definecolor{colorgray}{rgb}{0,0,0}

  \definecolor{coloroftheorems}{rgb}{0,0,0}
  \definecolor{colorigor}{rgb}{0,0,0}
  \definecolor{amethyst}{rgb}{0,0,0}
  
  \def\cole{\color{coloroftheorems}}
  
  \definecolor{colororange}{rgb}{0.8,0.2,0}
  \definecolor{colorpurple}{rgb}{0.6,0.0,0.6}
  
\fi

  \def\bea{\begin{align}}
  \def\ena{\end{align}}

\def\bega{\begin{aligned}}
  \def\enda{\end{aligned}}

\def\bcase{\begin{cases}}
  \def\ecase{\end{cases}}

\def\bmx{\begin{bmatrix}}
  \def\emx{\end{bmatrix}}

\def\hh{\textnormal{h}}
\def\cco{\textnormal{co}}
\def\ttan{\textnormal{tan}}

\begin{document}
\baselineskip=12.6pt

$\,$
\vskip1.2truecm
\title[A BKM-type criterion for the Euler equations]{A BKM-type criterion for the Euler equations}

\author[M.S.~Ayd\i n]{Mustafa Sencer Ayd\i n}
\address{Department of Mathematics, University of Southern California, Los Angeles, CA 90089}
\email{maydin@usc.edu}
%
%


\begin{abstract}
	We establish a new BKM-type blow-up criterion for solutions of the incompressible Euler equations that belong to Sobolev or H\" older spaces. Our criterion involves the $L^2$ norm in time of the $L^\infty$ norm of the first order tangential derivatives. Moreover, it applies to various domains such as the full space, the half-space, torus, (in)finite channel, and domains with curved boundaries. Additionally, we provide a mixed criterion involving the $L^1_t L^\infty(\Omega_1)$ norm of the vorticity and the $L^2_t L^\infty(\Omega_2)$ norm of the first order conormal derivatives of the velocity where
	$\Omega_1 \cup \Omega_2 = \Omega$ is a suitable decomposition of the physical space. Finally, we prove a blow-up criterion for the class of solutions that belong to the Sobolev conormal spaces that is recently constructed in~\cite{AK1}. 
\end{abstract}

\maketitle

\date{}

\startnewsection{Introduction}{sec01}

We consider the three-dimensional incompressible Euler equations
\begin{align}
	u_t + u\cdot \nabla u + \nabla p =0
	\comma \nabla \cdot u = 0
	\comma (x,t) \in \Omega \times (0,T)
	,
	\label{euler}
\end{align}
where $\Omega\subseteq\mathbb{R}^3$.
When $\partial \Omega \neq \emptyset$, the Euler equations are coupled with the slip boundary condition
\begin{align}
	u \cdot n = 0 
	\text{ on } \partial \Omega \times (0,T)
	,
	\label{eulerb}
\end{align}
where $n$ is the unit outward normal. 

In~\cite{BKM}, Beale, Kato, and Majda characterized the breakdown of strong solutions (in Sobolev spaces) of Euler equations involving the quantity 
\begin{align}
	\int_0^t \Vert \omega(s)\Vert_{L^{\infty}}\,ds
	,\label{EQ24}
\end{align} establishing that the solutions stay smooth as long as the $L^1_t L^\infty$ norm of the vorticity
stays finite. 
Next, Ponce in~\cite{P} proved a similar result using the norm on the deformation tensor
\begin{align}
	\sum_{i,j=1}^{3} \int_0^t \Vert \partial_i u_j(s)+ \partial_j u_i(s)\Vert_{L^{\infty}}\,ds
	.\llabel{EQ25}
\end{align}
Later, Constantin, Fefferman, and Majda utilized the direction of the vorticity 
to present a continuation (or blow-up) criterion involving
\begin{align}
	\int_0^t (1+\Vert u(s)\Vert_{L^{\infty}}) \left\Vert \nabla \frac{\omega(s)}{|\omega(s)|}\right\Vert_{L^{\infty}}\,ds
	.\llabel{EQ26}
\end{align}
The first blow-up criterion that included a weaker norm than $L^\infty$ was obtained by Kozono and Taniuchi in~\cite{KT} where the authors
showed that it suffices to control
\begin{align}
	\int_0^t \Vert \omega(s)\Vert_{BMO}\,ds
	\llabel{EQ27}
\end{align}
to prevent Euler solutions from blowing up. 
With Ogawa in~\cite{KOT}, they have also extended this result by establishing the criterion on
\begin{align}
	\int_0^t \Vert \omega(s)\Vert_{\dot{B}^0_{\infty,\infty}}\,ds
	,\llabel{EQ29}
\end{align}
which is indeed an improvement since $BMO \subset \dot{B}^0_{\infty,\infty}$.
Another result deviating from $L^\infty$ is due to Chae~\cite{Ch3}
who established that Euler solutions stay smooth as long as
\begin{align}
	\int_0^t \Vert \omega_\hh(s)\Vert_{\dot{B}^0_{\infty,1}}\,ds
	\llabel{EQ30}
\end{align}
is finite, where $\omega_\hh = (\omega_1,\omega_2)$. Although we have the embedding $\dot{B}^0_{\infty,1} \subset L^\infty$,
we note that this criterion only uses the horizontal components of the vorticity.
For the equations posed in $\mathbb{R}^3$, we also refer the reader to~\cite{Ch4, Ch5, CC2, CW2, K} for results on Sobolev solutions,~\cite{BD, CW3} for H\"older solutions,~\cite{Ch1, Ch2}
for Triebel-Lizorkin and Besov solutions,~\cite{CI, CK, CW1} for axisymmetric solutions,
and~\cite{Co, CC1, DHY1, DHY2, GHKR} for geometric criteria.
Finally, the BKM criterion for bounded domains is proven by Ferrari in~\cite{F}; see also~\cite{OT, SY, Z}.

In this work, we present a new BKM-type criterion involving the $L^2_t W^{1,\infty}_\ttan$ norm of the velocity (see Sections~\ref{sec02} and~\ref{sec04} for the definition of tangential Sobolev spaces and norms). 
All of the criteria above involve the normal derivative of either the vorticity, the direction of the vorticity, the deformation tensor, or the pressure. Therefore, the improvement here is that the criterion presented in Theorem~\ref{T03}
only involves tangential derivatives. Since we require the $L^2$-in-time norm, unlike the results we have mentioned that require the $L^1$-in-time norm, this improvement is at the expense of time integrability. It is not clear whether it is possible to improve our result and present a criterion involving $L^1_t W^{1,\infty}_\cco$ of the velocity.

We note that the BKM-type results involve an analysis of the vortex stretching term $\omega \cdot \nabla u$. 
This is also true in our case. We show that $\omega \cdot \nabla u$ behaves like a product of the vorticity with a tangential derivative of velocity plus a quadratic tangential derivative of $u$. In this way, the $L^2_t L^\infty$ control over the first
order tangential derivatives of $u$ inhibit vorticity accumulation by only allowing the vorticity to grow exponentially. Since we already know that the $L^1_t L^\infty$ norm of the vorticity governs the smoothness of solutions and we have a control on the $L^\infty_t L^\infty$ norm of the vorticity, we conclude that Euler solutions remain smooth as long as the velocity stays bounded in $L^2_t W^{1,\infty}_\cco$.

Our second main result is to glue the BKM criterion with the $L^2_t W^{1,\infty}_\cco$ norm of the velocity.
Under a suitable decomposition $\Omega= \Omega_1\cup\Omega_2$, we characterize the blow-up of solution
using the sum of $L^1_t W^{1,\infty}_\cco (\Omega_1)$ norm of the velocity 
(see Sections~\ref{sec02} and~\ref{sec04}) and the $L^1_t L^\infty (\Omega_2)$
norm of the vorticity. One may expect to achieve this using $W^{1,\infty}_\ttan$ instead of $W^{1,\infty}_\cco$ norm.
However, our proof shows that this may fail. In particular, when the closure of $\Omega_1$ intersects with the boundary
we end up with a conormal derivative of the form $\varphi\partial_n$,
where $\partial_n$ is the normal derivative and $\varphi$ is the distance to the boundary. 
Moreover, it is not possible to
control $\varphi\partial_n u$ locally using only the tangential derivatives due to the non-locality of the pressure term.
Therefore, in Theorems~\ref{T04} and~\ref{T05}, we consider $W^{1,\infty}_\cco$ rather than $W^{1,\infty}_\ttan$.

Finally, we provide a BKM-type criterion for the solutions of Euler equations belonging to Sobolev conormal spaces. The literature on the well-posedness of Euler equations in such spaces is brief. 
In~\cite{MR1}, Masmoudi and Rousset studied the inviscid problem with Navier-boundary conditions utilizing Sobolev conormal spaces. As a byproduct of their analysis, they have established the existence and uniqueness of Euler equations in Sobolev conormal spaces; see also~\cite{A, AK2}.  
Next, in~\cite{BILN}, the authors considered a weaker set of assumptions and showed the well-posedness of the Euler equations.
Finally, in~\cite{AK1} and~\cite{A}, we extended their result by removing any conormal differentiability assumptions on the normal derivative of the velocity in $L^2$ and $L^p$-based Sobolev spaces, respectively. 
A common ground for all of these results is that there is no second-order normal differentiability requirement, and this is only possible at the expense of having more tangential regularity than the Euler solutions belonging to the usual Sobolev spaces. Indeed, the classical well-posedness theory requires at least $2.5^+$ derivatives in all directions, whereas in~\cite{AK1}, we only need one in the normal direction but four in the tangential directions.

The criterion controlling the growth of conormal norms of the Euler solutions is given by the 
$L_t^1 W^{2,\infty}_\cco \cap L_t^2 W^{1,\infty}_\cco$ norm of the velocity. Therefore,
there is a discrepancy between the two criteria introduced for the Euler solutions in Sobolev and Sobolev conormal spaces
in the same way that there is a discrepancy in the number of normal and tangential derivatives
needed to construct solutions. 
Since we do not have a BKM-type criterion for the conormal solutions, we separately estimate each norm included in Theorem~\ref{T02}. An a~priori comparison of Theorem~\ref{T02} and Theorem~\ref{T03} suggests that the Euler solutions in conormal spaces are less likely to blow up.

The rest of our work is as follows. In the next section, we introduce the tangential and conormal spaces with norms 
when the underlying domain has flat or no boundary. Then, we present our main results. 
Next, in Section~\ref{sec03}, we prove Theorems~\ref{T03} and~\ref{T04}. 
Extensions of these result to general bounded domains are presented in Section~\ref{sec04}.
Finally, in Section~\ref{sec031}, we prove Theorem~\ref{T02} for the half-space.

\startnewsection{Preliminaries and main results}{sec02}

When $\Omega = \mathbb{R}^3_+$, we denote $x \in \Omega$ by $x=(x_1,x_2,z)= (x_\hh,z)$ and $(\partial_1,\partial_2)$ by $\nabla_\hh$. Next, conormal derivatives are given by 
\begin{align}
	(Z_1,Z_2,Z_3)=(\partial_1,\partial_2,\varphi(z)\partial_z)
	\comma \varphi(z)=\frac{z}{1+z}
	.\llabel{EQ31}
\end{align}
For $1\le p\le \infty$, we define 
the tangential Sobolev spaces by
\begin{align}
	W^{m,p}_\ttan=
	W^{m,p}_\ttan(\Omega)
	=
	\{f \in L^p(\Omega) : Z^\alpha f \in L^p(\Omega), \alpha \in \mathbb{N}_0^3, 0\le |\alpha|\le m, \alpha_3=0  \}
	,
	\label{tan.def}
\end{align} 
and the Sobolev conormal spaces by
\begin{align}
	W^{m,p}_\cco=
	W^{m,p}_\cco(\Omega)
	=
	\{f \in L^p(\Omega) : Z^\alpha f \in L^p(\Omega), \alpha \in \mathbb{N}_0^3, 0\le |\alpha|\le m  \}
	,
	\llabel{con.def}
\end{align} 
with the norms 
\begin{align}
	\begin{split}
		\Vert f\Vert_{W^{m,p}_\ttan(\Omega)}^p
		=&
		\sum_{|\alpha|\le m, \alpha_3=0} \Vert Z^\alpha f\Vert_{L^p(\Omega)}^p
		\comma 1\le p< \infty\\
		\\
		\Vert f\Vert_{W^{m,p}_\cco(\Omega)}^p
		=&
		\sum_{|\alpha|\le m} \Vert Z^\alpha f\Vert_{L^p(\Omega)}^p
		\comma 1\le p< \infty,
		\llabel{norm}
	\end{split}
\end{align}
and
\begin{align}
	\begin{split}
		\Vert f\Vert_{W^{m,\infty}_\ttan(\Omega)}
		=&
		\sum_{|\alpha|\le m, \alpha_3=0} \Vert Z^\alpha f\Vert_{L^p(\Omega)},
		\\
		\Vert f\Vert_{W^{m,\infty}_\cco(\Omega)}
		=&
		\sum_{|\alpha|\le m} \Vert Z^\alpha f\Vert_{L^\infty(\Omega)}.
		\llabel{norm2}
	\end{split}
\end{align}

Next, we consider the case $\Omega = \mathbb{R}^3$. In this case, conormal derivatives, spaces, and norms
are defined as above upon letting
\begin{equation}\llabel{EQ32}
	\varphi(z)=
	\begin{cases}
		\frac{z}{1+z}, & z\ge 0 \\
		\frac{z}{1-z}, & z < 0,
	\end{cases}
\end{equation}
Additionally, when the spatial domain is $\mathbb{T}^2\times [0,1]$ or $\mathbb{R}^2\times [0,1]$, we consider
\begin{equation}\llabel{EQ33}
	\varphi(z)=
	\min(z,1-z).
\end{equation}
 Finally, when $\Omega = \mathbb{T}^3 = [0,1]^3/\sim$, we set
\begin{equation}
	\varphi(z)=
	\sin(2\pi z).\llabel{EQ50}
\end{equation}
In both of these cases, we still consider $W^{k,p}_\ttan(\Omega)$ as defined in \eqref{tan.def}.
However, our results apply when we use derivatives in any fixed two directions in defining $W^{k,p}_\ttan(\Omega)$. 
The analogues of these spaces for curved bounded domains are defined in Section~\ref{sec04}.

We present a partial refinement of the Beale-Kato-Majda criterion for the strong solutions of the Euler equations.
Here, we consider the class of solutions $u$ such that
\begin{align}
	u \in C([0,T];H^s(\Omega))\cap C^1([0,T];H^{s-1}(\Omega)) \comma s > 2.5,
	\label{EQ07}
\end{align}
for some $T>0$.
Now, we state our first main result.
We note that the spatial domain $\Omega$ is either $\mathbb{R}^3$, $\mathbb{R}^3_+$,
$\mathbb{T}^3$, (in)finite channel, (in)finite cylinder, or an open bounded domain.

\cole
\begin{theorem}[Blow-up criterion for the solutions in Sobolev spaces]
	\label{T03}
	Let $s> 2.5$.
	If the maximal time of existence $T^*$ for the unique solution $u$ 
	of the Euler equations satisfying \eqref{EQ07} for $T<T*$ is finite, then 
	\begin{align}
		\int_0^{T*} \Vert u(s)\Vert_{W^{1,\infty}_\ttan}^2 \,ds = \infty. 
		\label{EQ.main1}
	\end{align} 
\end{theorem}
\colb

In particular, Theorem~\ref{T03} implies that the blow-up has to occur in tangential derivatives.

\begin{Remark}\label{R01}
	Our proof shows that we may extend our result to other classes of solutions for which 
	the BKM criterion \eqref{EQ24} is established. For example, the condition \eqref{EQ.main1} 
	also characterizes the blow-up in the class of H\"older solutions, $C^{k,\alpha}$ for
	$k \in \mathbb{N}$ and $\alpha \in (0,1)$. 
\end{Remark}

\begin{Remark}\label{R06}
	When $\partial \Omega = \emptyset$ or when $\partial \Omega$ is flat, we may replace $\Vert u(s)\Vert_{W^{1,\infty}_\ttan}^2$
	by $\Vert \nabla_\hh u(s)\Vert_{L^{\infty}}^2$ in \eqref{EQ.main1}; see Section~\ref{sec03}.
\end{Remark}

\begin{Remark}\label{R02}
	Due to the result of Elgindi in~\cite{E}, it is now known that there are Euler solutions in $C^{1,\alpha}$
	that develop singularities in finite time. We recall that the solution profile 
	in~\cite{E} is a vorticity in self-similar coordinates solving axisymmetric Euler equations with no swirl.
	Since the singularity formation of this profile does not distinguish directions, 
	and since the criterion \eqref{EQ.main1} involves tangential derivatives, our result is in line with
	the behavior of the solution profile in~\cite{E}; see also~\cite{CH}. 
\end{Remark}

In our next result, we show that the BKM criterion can be glued with a criterion involving $W^{1,\infty}_\cco$ norm of $u$.
First, we only consider the case where $\partial \Omega$ is flat or empty.
\begin{definition}\label{D01}
  Let $\Omega_1,\Omega_2 \subseteq \Omega$ be open and smooth and $\chi \in C^\infty(\Omega)$. We say that the triplet $(\Omega_1,\Omega_2,\chi)$
  is compatible when 
  	\begin{itemize}
  		\item[i.] $d(\Omega_1^\textnormal{c}, \Omega_2^\textnormal{c})>0$,
  		\item[ii.] $\chi = 1$ in $\Omega_2^\textnormal{c}$, $\chi = 0$ in 
  		$\Omega_1^\textnormal{c}$, and $\partial_1 \chi = \partial_2 \chi = 0$,     
  		\item[iii.] $d(\Omega_1,\{z=0\})>0$ when $\Omega = \mathbb{R}^3,\mathbb{T}^3$.
  	\end{itemize} 	
   
\end{definition}
For the half-space, an example of a such triplet is given by 
$\Omega_1 = \mathbb{R}^2 \times (0,b)$, $\Omega_2 = \mathbb{R}^2 \times (a,\infty)$,
where $0< a<b$, and $\chi = \chi(z)$ satisfying Definition~\ref{D01}(ii). When $\Omega = \mathbb{R}^3$, 
we may take $\Omega_1 = \mathbb{R}^2 \times \left((a,\infty) \cup (-\infty,-a)\right)$ and $\Omega_2 = \mathbb{R}^2 \times (-b,b)$.
Similar examples may be constructed for $\mathbb{T}^3$ or the (in)finite channel. The analogue of Definition~\ref{D01}
for bounded curved domains is given in Section~\ref{sec04}.

\cole
\begin{theorem}[A mixed blow-up criterion]
	\label{T04}
	Let $(\Omega_1,\Omega_2,\chi)$ be a compatible pair in the sense of Definition~\ref{D01}, and 
	assume that $s> 2.5$.
	If the maximal time of existence $T^*$ for the unique solution $u$ 
	of the Euler equations satisfying \eqref{EQ07} for $T<T*$ is finite, then 
	\begin{align}
		 \int_0^{T*} (\Vert u(s)\Vert_{W^{1,\infty}_\cco(\Omega_1)}^2+
		 \Vert \omega(s)\Vert_{L^{\infty}(\Omega_2)}) \,ds = \infty. 
		\llabel{EQ.main4}
	\end{align} 
\end{theorem}
\colb

See Theorem~\ref{T05} for the case when $\Omega$ is a general bounded domain. Next, we recall \cite[Theorem~2.1]{AK1} establishing the existence and uniqueness of the Euler equations in Sobolev conormal spaces.

\cole
\begin{theorem}(Euler equations in the Sobolev conormal spaces)
	\label{T01}
	Assume that $\Omega$ is $\mathbb{R}^3_+$, $\mathbb{T}^2\times [0,1]$, or $\mathbb{R}^2 \times [0,1]$. Let $u_0 \in H^4_\cco(\Omega) \cap W^{2,\infty}_\cco(\Omega) \cap W^{1,\infty}(\Omega)$
	be such that $\div u_0 = 0$, and $u_0 \cdot n= 0$ on~$\partial \Omega$.
	For some $T>0$,
	there exists  a unique solution 
	$u \in L^\infty(0,T;H^4_\cco(\Omega)\cap W^{2,\infty}_\cco(\Omega)\cap W^{1,\infty}(\Omega))$
	to the Euler equations
	\eqref{euler} and \eqref{eulerb}
	such that
	\begin{align}
		\sup_{[0,T]}
		(\Vert u(t)\Vert_{H^4_\cco}^2
		+\Vert u(t)\Vert_{W^{2,\infty}_\cco}^2
		+\Vert \nabla u(t)\Vert_{L^\infty}^2)
		\le M,
		\llabel{EQ.main3}
	\end{align}
	for $M>0$ depending on the norms of the initial data. 
\end{theorem}
\colb

As a consequence, we have the following regularity-in-time result.

\cole
\begin{Lemma}[Time continuity of solutions]\label{L03}
	Let $u_0 \in H^4_\cco(\Omega) \cap W^{2,\infty}_\cco(\Omega) \cap W^{1,\infty}(\Omega)$
	be such that $\div u_0 = 0$ in~$\Omega$, and $u_0 \cdot n= 0$ on~$\partial \Omega$.
	Then, the unique solution provided by Theorem~\ref{T01} belongs to~$C([0,T];H^4_\cco(\Omega)\cap W^{2,\infty}_\cco(\Omega) \cap W^{1,\infty}(\Omega)) \cap C^1([0,T];H^3_\cco(\Omega) \cap W^{1,\infty}_\cco(\Omega))$. 
\end{Lemma}
\colb

We now state our blow-up criterion for solutions belonging to the Sobolev conormal spaces.

\cole
\begin{theorem}[Blow-up criterion for solutions in Sobolev conormal spaces]
	\label{T02}
	Assume that $\Omega$ is $\mathbb{R}^3_+$, $\mathbb{T}^2\times [0,1]$, or $\mathbb{R}^2 \times [0,1]$.
	If the maximal time of existence $T^*$ for the unique solution $u$ 
	of the Euler equations given by Theorem~\ref{T01} and Lemma~\ref{L03} is finite, then 
	\begin{align}
		\int_0^{T*} \left(\Vert u(s)\Vert_{W^{1,\infty}_\cco}^2 + \Vert u(s)\Vert_{W^{2,\infty}_\cco}\right) \,ds = \infty. 
		\label{EQ.main2}
	\end{align} 
\end{theorem}
\colb

In particular, when the 
left-hand side of \eqref{EQ.main2} is finite, we may control the $H^4_\cco \cap W^{2,\infty}_\cco \cap W^{1,\infty}$
norm uniformly in time. 

\startnewsection{Proofs of Theorem~\ref{T03} and Theorem~\ref{T04} when $\partial \Omega$ is empty or flat}{sec03}

Our proofs rely on the fact that the uniform norm of the vorticity can only grow exponentially when $L^2_t W^{1,\infty}_\ttan$ norm of the velocity is bounded. 
First, we present the proof of Theorem~\ref{T03}.

\begin{proof}[proof of Theorem~\ref{T03}]
   It suffices to show that
    \begin{align}
   	\Vert \omega\Vert_{L^{\infty}}
   	\lec
   	\Vert \omega_0\Vert_{L^{\infty}}
   	+\int_0^t \left(\Vert \omega \Vert_{L^{\infty}} \Vert \nabla_\hh u\Vert_{L^\infty} + \Vert \nabla_\hh u\Vert_{L^\infty}^2\right) \,ds
   	,\label{EQ01}
   \end{align}
   for $t \in [0,T]$ and $T<T^*$.
   Indeed, suppose that \eqref{EQ01} holds and that 
   \begin{align}
   	 \int_0^{T*} \Vert u(s)\Vert_{W^{1,\infty}_\ttan}^2 \,ds \le C < \infty. 
   	\llabel{EQ100}
   \end{align} 
	Then, we may use \eqref{EQ01} and Gr\"onwall's inequality to conclude that $\Vert \omega\Vert_{L^{\infty}}$
	stays bounded uniformly up to the time $T^*$. Finally, since the $L^1(0,T^*;L^\infty)$ norm of the vorticity is 
	bounded by the $L^\infty(0,T^*;L^\infty)$-norm, Theorem~\ref{T03} follows from the BKM criterion.

  To prove \eqref{EQ01}, we perform $L^p$ estimates on the vorticity equation obtaining
  \begin{align}
  	\frac{1}{p}\frac{d}{dt}\Vert \omega\Vert_{L^{p}}^p
  	=\int_\Omega \omega \cdot \nabla u \omega |\omega|^{p-2}
  	\le \Vert u \cdot \nabla \omega\Vert_{L^{p}}\Vert \omega\Vert_{L^{p}}^{p-1}
  	,\llabel{EQ90} 
  \end{align}
  where we have used that 
  \begin{align}
  	\int_\Omega u\cdot \nabla \omega \omega |\omega|^{p-2}
  	 =\frac{1}{p}\int_\Omega \nabla \cdot (u |\omega|^p)
  	 =0.\llabel{EQ91} 
  \end{align} 
 Dividing by $\Vert \omega\Vert_{L^{p}}^{p-1}$ and sending $p\to \infty$, we arrive at
 \begin{align}
 	\Vert \omega\Vert_{L^{\infty}}
 	 \lec
 	  1+\int_0^T \Vert \omega \cdot \nabla u\Vert_{L^{\infty}} \,dt
 	  ,\llabel{EQ02}
 \end{align} 
  for $T<T^*$.
  Since 
  \begin{align}
  	\Vert \omega_3\Vert_{L^{\infty}} \lec \Vert \nabla_\hh u\Vert_{L^\infty}
  	\llabel{EQ03}
  \end{align}
  and 
  \begin{align}
  	\Vert \partial_z u\Vert_{L^{\infty}} \lec \Vert \omega\Vert_{L^{\infty}}+\Vert \nabla_\hh u\Vert_{L^\infty}
  	,\label{EQ04}
  \end{align}
  we may estimate
  \begin{align}
  	\Vert \omega_3 \cdot \partial_z u\Vert_{L^{\infty}}
  	\lec
  	\Vert \nabla_\hh u\Vert_{L^\infty}\Vert \omega\Vert_{L^{\infty}}+\Vert \nabla_\hh u\Vert_{L^\infty}^2
  	\label{EQ05}
  \end{align}
  while
  \begin{align}
  	\Vert \omega_\hh \cdot \nabla_\hh u\Vert_{L^{\infty}}
  	\lec
  	\Vert \omega\Vert_{L^{\infty}}\Vert \nabla_\hh u\Vert_{L^\infty}
  	,\label{EQ06}
  \end{align}
  from where \eqref{EQ01} follows.   	
  \end{proof}

Next, we prove Theorem~\ref{T04}.

\begin{proof}[proof of Theorem~\ref{T04}]
	Our aim is to control the $L^1_t L^\infty$-norm of the vorticity.
	We start by writing
	\begin{align}
		\partial_t (\chi \omega)
		 +u\cdot \nabla (\chi \omega)
		  = \chi \omega \cdot \nabla u
		   + u \cdot \nabla \chi \omega
		   ,\llabel{EQ80}
	\end{align}
	from where we obtain
	\begin{align}
		\Vert \chi \omega\Vert_{L^{\infty}(\Omega)}
		 \lec
		  1+\int_0^T (\Vert \chi \omega \cdot \nabla u\Vert_{L^{\infty}(\Omega_1)}
		   +\Vert u \cdot \nabla \chi \omega\Vert_{L^{\infty}(\Omega_1)})\,ds
		   .\label{EQ81}
	\end{align}
Keeping $\chi$ attached, we may estimate the vortex stretching term as
\begin{align}
	\Vert \chi \omega \cdot \nabla u\Vert_{L^{\infty}(\Omega_1)}
	 \lec
	  \Vert \chi \omega\Vert_{L^{\infty}(\Omega_1)}
	   \Vert \nabla_\hh u\Vert_{L^{\infty}(\Omega_1)}
	    +\Vert \nabla_\hh u\Vert_{L^{\infty}(\Omega_1)}^2
	    ,\label{EQ82}
\end{align}
where we have utilized \eqref{EQ05} and \eqref{EQ06} upon replacing $\omega$ by $\chi \omega$. 
Since the tangential derivatives of $\chi$ vanish, we obtain
\begin{align}
	\Vert u \cdot \nabla \chi \omega\Vert_{L^{\infty}(\Omega_1)}
	 =\Vert u_3 \partial_z \chi \omega\Vert_{L^{\infty}(\Omega_1)}
	  \lec
	   \Vert u_3 \partial_z u\Vert_{L^{\infty}(\Omega_1)}
	    +\Vert u_3 \nabla_\hh u\Vert_{L^{\infty}(\Omega_1)}
	     \lec
	      \Vert u\Vert_{W^{1,\infty}_\cco(\Omega_1)}^2
	    .\label{EQ83}
\end{align}
In the last inequality,  
we have utilized the slip boundary condition, Hardy's inequality, and the
incompressibility by writing
\begin{align}
	\Vert u_3 \partial_z u\Vert_{L^{\infty}(\Omega_1)}
	 =\left\Vert \frac{u_3}{\varphi(z)} \varphi(z)\partial_z u\right\Vert_{L^{\infty}(\Omega_1)}
	 \lec
	  \Vert \partial_z u_3\Vert_{L^{\infty}(\Omega_1)}\Vert Z_3 u\Vert_{L^{\infty}(\Omega_1)}
	   \lec 
	    \Vert u\Vert_{W^{1,\infty}_\cco(\Omega_1)}^2
	    ,\llabel{EQ84}
\end{align}
when $\partial \Omega \neq \emptyset$.
On the other hand, for the case where $\Omega = \mathbb{R}^3, \mathbb{T}^3$
we employ Definition~\ref{D01}(iii) and obtain 
$\Vert u_3/\varphi\Vert_{L^{\infty}(\Omega_1)}\lec \Vert u_3\Vert_{L^{\infty}(\Omega_1)}$.
Now, combining \eqref{EQ81}--\eqref{EQ83}, we have
	\begin{align}
	\Vert \chi \omega\Vert_{L^{\infty}(\Omega)}
	\lec
	1+\int_0^T (\Vert \chi \omega \Vert_{L^{\infty}(\Omega_1)}\Vert u \Vert_{W^{1,\infty}_\cco(\Omega_1)}
	+\Vert u \Vert_{W^{1,\infty}_\cco(\Omega_1)}^2)\,ds
	,\label{EQ85}
\end{align}
for $T<T^*$.
Therefore, when we assume that
\begin{align}
	 \int_0^{T*} (\Vert u\Vert_{W^{1,\infty}_\cco(\Omega_1)}^2+
	\Vert \omega\Vert_{L^{\infty}(\Omega_2)}) \,ds < \infty,
	\label{EQ86}
\end{align} 
Gr\"onwall's inequality and \eqref{EQ85} imply that
\begin{align}
	\Vert \chi \omega\Vert_{L^{\infty}(\Omega_1)} \lec 1
	,\label{EQ87}
\end{align}
where the implicit constant may depend on $T^*$. 
Finally, we note that for $T<T^*$
\begin{align}
	\int_0^T \Vert \omega\Vert_{L^{\infty}(\Omega)} \,dt
	 \lec
	  \int_0^T \Vert \chi\omega\Vert_{L^{\infty}(\Omega)} \,dt
	   +\int_0^T \Vert (1-\chi)\omega\Vert_{L^{\infty}(\Omega)} \,dt
	   \lec
	    1+\int_0^T \Vert (1-\chi)\omega\Vert_{L^{\infty}(\Omega_2)} \,dt
	     \lec 1
	     ,\llabel{EQ88}
\end{align}
where we have used \eqref{EQ87} in the second inequality and \eqref{EQ86}
in the third. Letting $T\to T^*$, we obtain a contradiction to the BKM criterion.
\end{proof}

\startnewsection{Extension to domains with curved boundaries}{sec04}

Let $\Omega\subseteq \mathbb{R}^3$ be a smooth, open, and bounded domain.
We fix a finite open cover $\{U_i\}_i$ of $\Omega$,  
\begin{align}
	\Omega = \bigcup_{i=0}^{m_0} U_i
	,\llabel{EQ63}  
\end{align}
where $U_i \cap \partial \Omega = \emptyset$ if and only if $i=0$.
By utilizing a partition of unity $\{\phi_i\}_i$ subordinate to $\{U_i\}_i$, we drop the index $i$
and consider $(U,\phi)=(U_i,\phi_i)$ for $i\ge1$. Then, upon rotation and translation, 
$U \cap \partial \Omega$ is given by a graph of a smooth function $\psi$. 
Next, we fix an orthonormal frame $(\tau,\bar{\tau},n)$ where $n$ is the unit normal vector and write
\begin{align}
	(\partial_1^\psi,\partial_2^\psi,\partial_3^\psi)
	=
	(\partial_\tau, \partial_{\bar{\tau}}, \partial_n)
	=
	(\tau \cdot \nabla, \bar{\tau}\cdot \nabla, n \cdot \nabla)
,\llabel{EQ64}
\end{align}
and \begin{align}
	(u_1^\psi,u_2^\psi,u_3^\psi)
	=
	(u_\tau, u_{\bar{\tau}}, u_n)
	=
	(\tau \cdot u, \bar{\tau}\cdot u, n \cdot u)
	.\llabel{EQ38}
\end{align}
Then, we also introduce $(g_{ij})_{3\times 3}$ such that 
\begin{align}
	\partial_i = g_{ij} \partial_j^\psi
	,\text{ and }  u_i = g_{ij} u_j^\psi\comma i=1,2,3.
	\llabel{EQ33}
\end{align}
Additionally, we have $g^{-1}=g^T$.
Tangential and conormal derivatives are given by 
$(\partial_1^\psi,\partial_2^\psi)$ and 
\begin{align}
	(Z_1,Z_2,Z_3)=(\partial_1^\psi,\partial_2^\psi,\varphi(x)\partial_3^\psi)
	,\llabel{EQ34}
\end{align}
respectively,
where $\varphi:\Omega\to \mathbb{R}_+$ is the distance to the boundary.
When $(U,\phi)= (U_0,\phi_0)$, we consider 
$Z_i=\partial_i$ for $i=1,2,3.$

For $1\le p\le \infty$,
we define the tangential Sobolev spaces 
to be the set of functions $f \in L^p(\Omega)$ such that
\begin{align}
	\Vert f\Vert_{W^{k,p}_\ttan}
	 =
	  \Vert f\Vert_{W^{k,p}(U_0)}+
	   \sum_{i=1}^{m_0} \sum_{|\alpha|\le k, \alpha_3=0} \Vert Z^\alpha f\Vert_{L^{p}(U_i)}
	    < \infty
	   .\llabel{EQ65}
\end{align}
Similarly, we define the Sobolev conormal spaces to be the set of $f\in L^p(\Omega)$
such that
\begin{align}
	\Vert f\Vert_{W^{k,p}_\cco}
	=
	\Vert f\Vert_{W^{k,p}(U_0)}+
	\sum_{i=1}^{m_0} \sum_{|\alpha|\le k} \Vert Z^\alpha f\Vert_{L^{p}(U_i)}
	< \infty
	.\llabel{EQ66}
\end{align}

Next, we show that \eqref{EQ04} is still true on curved domains. 

\cole
\begin{Lemma}[Estimates on $\partial_n u$]\label{L02}
	Let $U \in \{U_i\}_i$ with $U \neq U_0$. Assume that $u$ is
	a smooth divergence-free vector field on $\Omega$ and $\omega = \curl u$. Then, 
	\begin{align}
		\Vert \partial_n u\Vert_{L^{\infty}(U)}
		\lec
		\Vert \omega\Vert_{L^{\infty}(U)}+\Vert u\Vert_{W^{1,\infty}_\ttan(U)}
		.\label{EQ39}
	\end{align}
\end{Lemma}
\colb

We note that this claim was proven by~\cite[Lemma~7]{BILN}
and by~\cite{MR1} in slightly different forms. In addition, 
when $U=U_0$, the inequality \eqref{EQ39} may be obtained as in Section~\ref{sec03}.

\begin{proof}[proof of Lemma~\ref{L02}]
	We first note that the divergence-free condition in $U$ reads
	\begin{align}
		0 = \partial_i u_i
		= g_{ij}\partial_j^\psi (g_{ik} u_k^\psi)
		=g_{ij}\partial_j^\psi g_{ik} u_k^\psi
		+g_{ij} g_{ik}\partial_j^\psi  u_k^\psi
		,\llabel{EQ41}
	\end{align}
	from where using that $g_{ij} g_{ik} = \delta_{ij}$, we arrive at
	\begin{align}
		\partial_n u_n =\partial_3^\psi  u_3^\psi 
		=-g_{ij}\partial_j^\psi g_{ik} u_k^\psi
		+\partial_1^\psi  u_1^\psi+\partial_2^\psi  u_2^\psi
		.\label{EQ42}
	\end{align}
	Therefore, the normal derivative of the normal component of $u$ consists of $u$ and
	its tangential derivatives. Next, we expand $\omega$ as 
	\begin{align}
		\omega_i 
		=
		\epsilon_{ijk} \partial_j u_k
		=
		\epsilon_{ijk} g_{jl} \partial_l^\psi (g_{km} u_m^\psi)
		=
		\epsilon_{ijk} g_{j3} g_{k1} \partial_3^\psi u_1^\psi
		+\epsilon_{ijk} g_{j3} g_{k2} \partial_3^\psi u_2^\psi
		+ L_i
		,\llabel{EQ43}
	\end{align}
	where $L_i$ only contains $u$ and its tangential derivatives, i.e.,
	\begin{align}
		\sum_{i=1}^{3}\Vert L_i\Vert_{L^{\infty}(U)} \lec \Vert u\Vert_{W^{1,\infty}_\ttan(U)}
		.\label{EQ44}
	\end{align}
	Now, we rewrite $\omega$ as
	\begin{align}
		\omega - L = g_{\cdot 3}\times g_{\cdot 1} \partial_3^\psi u_1^\psi 
		+g_{\cdot 3}\times g_{\cdot 2} \partial_3^\psi u_2^\psi
		,\label{EQ45}
	\end{align}
	where $g_{\cdot i}$ is the $i^\text{th}$ column of $g$.
	Finally, taking the dot product of \eqref{EQ45} with 
	$g_{\cdot 2}$ and $g_{\cdot 1}$, respectively, we obtain
	\begin{align}
		\partial_3^\psi u_1 = (\omega - L) \times g_{\cdot 2} 
		,\text{ and } 
		\partial_3^\psi u_2 = (\omega - L) \times g_{\cdot 1}
		,\label{EQ46} 
	\end{align}
	from where, we conclude \eqref{EQ39}. 
\end{proof}

We are in a position to prove Theorem~\ref{T03} for bounded domains.

\begin{proof}[proof of Theorem~\ref{T03}]
	We begin by expanding the vortex stretching term in $U$ by writing
\begin{align}
	\omega_i \partial_i u
	=
	\epsilon_{ijk} \partial_j u_k \partial_i u
	=
	\epsilon_{ijk} g_{jl} \partial_l^\psi (g_{km} u_m^\psi) g_{ii'} \partial_{i'}^\psi u
	.\llabel{EQ35}
\end{align}
We first note that when $i' = 3 = l$, we arrive at a term that
involves 
\begin{align}
	\epsilon_{ijk} g_{j3} g_{i3}
	= \epsilon_{ijk} \delta_{ij}
	=0
	.\llabel{EQ47}
\end{align}
Therefore, we have
\begin{align}
	\omega_i \partial_i u
	= 
	\sum_{\alpha=1}^{2}\epsilon_{ijk} g_{j3} g_{km} \partial_3^\psi u_m^\psi g_{i\alpha} \partial_{\alpha}^\psi u
	+
	\sum_{\beta=1}^{2}\epsilon_{ijk} g_{j\beta} \partial_\beta^\psi (g_{km} u_m^\psi) g_{i3} \partial_{3}^\psi u
	+\bar{L}
	,\label{EQ48}
\end{align}
where $\Vert \bar{L}\Vert_{L^{\infty}(U)} \lec \Vert u\Vert_{W^{1,\infty}_\ttan(U)}^2$.
Finally, employing Lemma~\ref{L02}, it follows from \eqref{EQ48} that
\begin{align}
	\Vert \omega\cdot \nabla u\Vert_{L^{\infty}(U)}
	\lec
	\Vert \omega\Vert_{L^{\infty}(U)}\Vert u\Vert_{W_\ttan^{1,\infty}(U)}
	+\Vert u\Vert_{W^{1,\infty}_\ttan(U)}^2
	.\label{EQ49}
\end{align}

Now, recalling that $\{\phi_i\}_i$ is a partition of unity subordinate to $\{U_i\}_i$,
we arrive at
\begin{align}
	\partial_t(\phi_i \omega) + u\cdot \nabla (\phi_i\omega) = \phi_i \omega \cdot \nabla u + u\cdot \nabla \phi_i \omega
	.\label{EQ66}
\end{align}
Next, we perform $L^\infty$ estimates up to time $T<T^*$ obtaining
\begin{align}
	\begin{split}
	\Vert \phi_i \omega\Vert_{L^{\infty}(\Omega)}
	 &\lec
	  1+\int_0^T \Vert \phi_i \omega\cdot \nabla u\Vert_{L^{\infty}(U_i)}
	  +\Vert u\cdot \nabla \phi_i \omega\Vert_{L^{\infty}(U_i)}\,dt
	  \\&\lec
	  1+\int_0^T (\Vert \omega\Vert_{L^{\infty}(U_i)}\Vert u\Vert_{W_\ttan^{1,\infty}(U_i)}
	  +\Vert u\Vert_{W^{1,\infty}_\ttan(U_i)}^2
	  )\,dt
	  ,\llabel{EQ67}
	  \end{split}
\end{align}
where we have used \eqref{EQ49} and the boundedness of $\phi$ in the second inequality.
Summing over $i$ and using that $\sum_i \phi_i = 1$, we arrive at
\begin{align}
		\Vert \omega\Vert_{L^{\infty}(\Omega)}
		\lec
		1+\int_0^T (\Vert \omega\Vert_{L^{\infty}(\Omega)}\Vert u\Vert_{W_\ttan^{1,\infty}(\Omega)}
		+\Vert u\Vert_{W^{1,\infty}_\ttan(\Omega)}^2
		)\,dt
		,\label{EQ68}
\end{align}
for all $T<T^*$.
With \eqref{EQ68} established, 
we may utilize Gr\"onwall's inequality and conclude the proof.
\end{proof}

Next, we proceed to establishing an analogue of Theorem~\ref{T04} for bounded domains.

\begin{definition}\label{D02}
	Let $\Omega\subset \mathbb{R}^3$ be a smooth, open and bounded domain. Then, $\Omega$
	admits a compatible covering if there exists a finite open covering $\{U_i\}_i$ of $\Omega$
	and a partition of unity $\{\phi_i\}_i$ subordinate to $\{U_i\}_i$ such that
	each $\phi_i$ has vanishing tangential derivatives. 
\end{definition} 

Examples of such domains include open balls, periodic cylinders and their diffeomorphic copies.
Indeed, for $B=B(0,1)$ we may take $U_1 = B\setminus B(0,a)$ and $U_2 = B(0,b)$, where $0<a<b$,
with radial and smooth cutoff functions. Now, we state and prove our result on the mixed BKM-type
criterion for bounded domains.

\cole
\begin{theorem}[A mixed blow-up criterion for bounded domains]
	\label{T05}
	Let $\Omega \subseteq \mathbb{R}^3$ be an open, smooth, and bounded domain admitting a compatible covering 
	$(\{U_i\}_i,\{\phi_i\}_i)$ in the sense of Definition~\ref{D02}. Moreover,
	for a partition $(\mathcal{U}_1,\mathcal{U}_2)$ of the collection $\{U_i\}_i$ let
	\begin{align}
		\mathcal{O}_i = \bigcup_{U \in \mathcal{U}_i} U \comma i=1,2
		.\llabel{EQ69}
	\end{align}
	If the maximal time of existence $T^*$ for the unique solution $u$ 
	of the Euler equations satisfying \eqref{EQ07} for $T<T*$ is finite, then 
	\begin{align}
		\int_0^{T*} (\Vert u(s)\Vert_{W^{1,\infty}_\cco(\mathcal{O}_1)}^2+
		\Vert \omega(s)\Vert_{L^{\infty}(\mathcal{O}_2)}) \,ds = \infty. 
		\label{EQ.main6}
	\end{align} 
\end{theorem}
\colb

\begin{proof}[proof of Theorem~\ref{T05}]
	We assume that the left-hand side of \eqref{EQ.main6} is finite. 
	Then, we bound the $L^1_T L^\infty(\Omega)$ norm of $\omega$ for $T<T^*$.
	Finally, we conclude by employing the BKM criterion.
	
	We let $U_i \in \mathcal{U}_1$ and multiply the Euler equations by $\phi_i$ obtaining
	\eqref{EQ66}. Furthermore, we may write
   \begin{align}
   	\begin{split}
   		\Vert \phi_i \omega\Vert_{L^{\infty}(U_i)}
   		&\lec
   		1+\int_0^T (\Vert \phi_i \omega\cdot \nabla u\Vert_{L^{\infty}(U_i)}
   		+\Vert u \cdot \nabla \phi_i \omega\Vert_{L^{\infty}(U_i)})\,dt
   		.\label{EQ70}
   	\end{split}
   \end{align}
   	Recalling \eqref{EQ45}, we note that
   	\begin{align}
   		u \cdot \nabla \phi_i \omega
   		  =  u \cdot \nabla \phi_i (g_{\cdot 3}\times g_{\cdot 1} \partial_3^\psi u_1^\psi 
   		  +g_{\cdot 3}\times g_{\cdot 2} \partial_3^\psi u_2^\psi
   		  ) + u \cdot \nabla \phi_i L
   		  .\label{EQ71}
   	\end{align}
   	We utilize that $\partial_j^\psi \phi_i = 0$ for $j=1,2$ and obtain
   	\begin{align}
   		u\cdot \nabla \phi_i
   		 =u_j \partial_j \phi_i
   		  =g_{jk} u^\psi_k g_{jl} \partial_l^\psi \phi_i
   		   =g_{jk} u^\psi_k g_{j3} \partial_3^\psi \phi_i
   		   =\delta_{k3} u^\psi_k \partial_3^\psi \phi_i
   		   =u^\psi_3 \partial_3^\psi \phi_i.\llabel{EQ72}
   	\end{align}
   	Combining this with \eqref{EQ71} and inserting the distance function $\varphi$, we arrive at
   	\begin{align}
   		u \cdot \nabla \phi_i \omega
   		=  \frac{u^\psi_3}{\varphi} \partial_3^\psi \phi_i (g_{\cdot 3}\times g_{\cdot 1} Z_3 u_1^\psi 
   		+g_{\cdot 3}\times g_{\cdot 2} Z_3 u_2^\psi
   		) + u \cdot \nabla \phi_i L
   		.\llabel{EQ73}
   	\end{align}
   	To estimate this term, we employ \eqref{EQ42}, Hardy's inequality, and \eqref{EQ44}
   	from where it follows that
   	\begin{align}
   		\Vert u \cdot \nabla \phi_i \omega\Vert_{L^{\infty}(U_i)}
   		 \lec
   		  \Vert u\Vert_{W^{1,\infty}_\cco}^2
   		  .\label{EQ74} 
   	\end{align}
   	It remains to estimate the vortex strecthing term for which we write
   	\begin{align}
   		\Vert \phi_i \omega\cdot \nabla u\Vert_{L^{\infty}(U_i)}
   		 \lec
   		  \Vert u\Vert_{W^{1,\infty}_\ttan(U_i)}
   		   +\sum_{\alpha =1}^{2} \Vert \phi_i \partial_3^\psi u \partial_\alpha^\psi u\Vert_{L^{\infty}(U_i)}
   		   ,\llabel{EQ75}
   	\end{align}
   	where we have used \eqref{EQ48}. Now, we employ \eqref{EQ46} for the last term
   	in the above inequality and obtain
   	\begin{align}
   		\Vert \phi_i \omega\cdot \nabla u\Vert_{L^{\infty}(U_i)}
   		\lec
   		\Vert u\Vert_{W^{1,\infty}_\ttan(U_i)}^2
   		+\Vert \phi_i \omega\Vert_{L^{\infty}(U_i)}\Vert u\Vert_{W^{1,\infty}_\ttan(U_i)}.\label{EQ76}
   	\end{align}
   	Combining \eqref{EQ70}, \eqref{EQ74}, and \eqref{EQ76}, we conclude
   	\begin{align}
   		\begin{split}
   			\Vert \phi_i \omega\Vert_{L^{\infty}(U_i)}
   			&\lec
   			1+\int_0^T (\Vert \phi_i \omega\Vert_{L^{\infty}(U_i)}\Vert u \Vert_{W^{1,\infty}_\cco(U_i)}
   			+\Vert u \Vert_{W^{1,\infty}_\cco(U_i)}^2)\,dt
   			,\label{EQ77}
   		\end{split}
   	\end{align}
   	for $T < T^*$.
   	
   	When the left-hand side of \eqref{EQ.main6} is finite, we may sum \eqref{EQ77} over 
   	$i$ satisfying $U_i \in \mathcal{U}_1$, and utilize Gr\"onwall's inequality and obtain
   	\begin{align}
   		\limsup_{T \to T^*} \sup_{[0,T]}\sum_{U_i \in \mathcal{U}_1}\Vert \psi_i \omega\Vert_{L^{\infty}(\Omega)}\lec 1
   		,\label{EQ78}
   	\end{align}
   	where the implicit constant may depend on $T^*$. Therefore, we may estimate the $L^1 L^\infty$ norm of $\omega$ as
   	\begin{align}
   		\int_0^T \Vert \omega\Vert_{L^{\infty}(\Omega)}\,dt
   		 \lec
   		  \sum_{j=1}^{2}\sum_{U_i \in \mathcal{U}_j} \int_0^T \Vert \phi_i \omega\Vert_{L^{\infty}(\Omega)}\,dt
   		   \lec
   		    1+\int_0^T \Vert \omega\Vert_{L^{\infty}(\mathcal{O}_2)}\,dt
   		     \lec 1
   		      ,\llabel{EQ79}
   	\end{align}
   	where we have used \eqref{EQ78} in the second inequality and the finiteness of \eqref{EQ.main6}
   	in the last. Upon letting $T \to T^*$, we may conclude the proof.   	
\end{proof}

\startnewsection{Proof of Theorem~\ref{T02}}{sec031}

Throughout this section, we assume that $\Omega = \mathbb{R}^3_+$. 
First, we note that
Lemma~\ref{L03} follows from our construction in~\cite{AK1} under standard considerations.
We may establish the weak continuity of the solution upon utilizing the Gel'fand
triple, $X \subseteq L^2 \subseteq X'$, where $X$
denotes $H^4_\cco(\Omega) \cap W^{2,\infty}_\cco(\Omega) \cap W^{1,\infty}(\Omega)$.
Here, we employ the convergence of the approximate solutions along with the uniform bounds we have obtained.
Since the Euler equations
are time-reversible and invariant under time shifts,
it suffices to prove the convergence in the $X$-norm as $t \to 0^+$ to obtain the continuity in $X$, and this convergence is a consequence of our estimates in~\cite{AK1}. 
Finally, for $t\in [0,T]$, we may show that $\partial_t u \in H^3_\cco \cap W^{1,\infty}_\cco$ upon conormalizing the normal derivative in the advective term. 
Now, we present the proof of Theorem~\ref{T02}.

\begin{proof}[proof of Theorem~\ref{T02}]
	We establish Theorem~\ref{T03} utilizing a contradiction argument. 
	We assume that $T^*<\infty$, and that 
	\begin{align}
		\int_0^{T*} \left(\Vert u(s)\Vert_{W^{1,\infty}_\cco}^2 
		 + \Vert u(s)\Vert_{W^{2,\infty}_\cco}\right) \,ds < \infty. 
		\label{EQas}
	\end{align}
	Then, we show that the solution $u$ stays bounded in $H^4_\cco \cap W^{2,\infty}_\cco \cap W^{1,\infty}$
	uniformly until the blow-up time $T^* < \infty$. Finally, by employing a continuation argument, we contradict that
	$T^*$ is maximal.
	
	For the rest of the proof, we only present the estimates on the a~priori level which we may justfiy by approximating the solution as in~\cite{AK1}.  
	Now, \eqref{EQ01} and Gr\"onwall's inequality implies
	\begin{align}
			\Vert \omega\Vert_{L^{\infty}}
			 \lec
			  \left(\Vert \omega_0\Vert_{L^{\infty}}+\int_0^T \Vert u\Vert_{W^{1,\infty}_\cco}^2\,ds\right)
			   \exp\left(T\int_0^T \Vert u\Vert_{W^{1,\infty}_\cco}^2 \,ds\right)
			   ,\llabel{EQ60}
	\end{align}
	for all $T<T^*$. Next, we use \eqref{EQas} and write
	\begin{align}
		\limsup_{T \to T^*}\sup_{[0,T]}\Vert \omega\Vert_{L^{\infty}} \lec 1
		,\llabel{EQ61}
	\end{align}
	where the implicit constant here and the rest of this proof may depend on $T^*$.
	In addition, we have that
	\begin{align}
		\Vert \nabla u(t)\Vert_{L^{\infty}} \le \Vert \omega(t)\Vert_{L^{\infty}}+\Vert u(t)\Vert_{W^{1,\infty}_\cco}
		 \lec 1+\Vert u(t)\Vert_{W^{1,\infty}_\cco}
		 ,\label{EQ09}
	\end{align}
	for $t \in [0,T]$ and $T < T^*$.
	Now, the pressure estimate given in~\cite[Proposition~3.2]{AK1} reads
	\begin{align}
		\Vert D^2 p(t)\Vert_{H^3_\cco}+\Vert \nabla p(t)\Vert_{H^3_\cco}
		\lec 
		\Vert u(t)\Vert_{H^4_\cco}(\Vert u(t)\Vert_{W^{2,\infty}_\cco}+\Vert \nabla u(t)\Vert_{L^\infty}+1)
		,
		\label{EQ.Pre}
	\end{align}
	for $t\in [0,T]$ and $T<T^*$.
	Combining this with the conormal derivative estimates in~\cite[Proposition~3.1]{AK1}, we write
	\begin{align}
		\sup_{[0,T]}\Vert u\Vert_{H^4_\cco}^2
		\lec
		\Vert u_0\Vert_{H^4_\cco}^2
		+\int_0^T
		\Bigl(
		\Vert u(s)\Vert_{H^4_\cco}^2 (
		\Vert u(s)\Vert_{W^{1,\infty}}+
		\Vert u(s)\Vert_{W^{2,\infty}_\cco}+1)		
		\Bigr)\,ds
		,
		\llabel{EQ.Con}
	\end{align}
	for $T<T^*$.
	Furthermore, the bound on $\nabla u$ in \eqref{EQ09} implies
	\begin{align}
		\sup_{[0,T]}\Vert u\Vert_{H^4_\cco}^2
		\lec
		1+\int_0^T
		\Bigl(
		\Vert u(s)\Vert_{H^4_\cco}^2 (
		\Vert u(s)\Vert_{W^{2,\infty}_\cco}+1)		
		\Bigr)\,ds
		,
		\llabel{EQ09}
	\end{align}
	for $T<T^*$.
	Finally, recalling \eqref{EQas}, we obtain that
	\begin{align}
		\limsup_{T \to T^*}\sup_{[0,T]}\Vert u\Vert_{H^4_\cco}^2
		\lec
		1,
		\label{EQ10}
	\end{align}
	upon utilizing the Gr\"onwall's inequality.
	
	Next, we establish that $\nabla u$ stays bounded in $L^\infty$ until time $T^*$.
	To achieve this, we consider a conormal derivative $Z$ and apply it to the Euler equations and write
	\begin{align}
		Zu_t + u\cdot \nabla Zu = -Zu \cdot \nabla u - (u \cdot Z \nabla u - u\cdot \nabla Zu) - Z \nabla p
		,\llabel{EQ11}
	\end{align}
	from where we may obtain
	\begin{align}
		\Vert Zu(t)\Vert_{L^{\infty}}
		 \lec
		  \Vert Zu_0\Vert_{L^{\infty}}
		  +\int_0^t (\Vert Zu \cdot \nabla u\Vert_{L^{\infty}}
		   +\Vert u \cdot Z \nabla u - u\cdot \nabla Zu\Vert_{L^{\infty}}
		    +\Vert Z \nabla p\Vert_{L^{\infty}}) \,ds
		    ,\label{EQ12}
	\end{align}
	for $t \in [0,T]$ and $T<T^*$.
    We bound the first two terms in the integral by the estimate
    \begin{align}
    	\Vert Zu \cdot \nabla u\Vert_{L^{\infty}}
    	+\Vert u \cdot Z \nabla u - u\cdot \nabla Zu\Vert_{L^{\infty}}
    	 \lec 
    	  \Vert \nabla u\Vert_{L^{\infty}}\Vert u\Vert_{W^{1,\infty}_\cco}
    	  ,\label{EQ13}
    \end{align}
while for the pressure term, we recall the inequality (see~\cite{MR2})
\begin{align}
	\Vert f\Vert_{L^{\infty}}\lec
	 \Vert \partial_z f\Vert_{H_{\textnormal{tan}}^{s_1}}^\frac12 \Vert f\Vert_{H_{\textnormal{tan}}^{s_2}}^\frac12
	 \comma s_1+s_2 >2,\llabel{EQ14}
\end{align}
where $H^s_{\textnormal{tan}}$ stands for the tangential Sobolev space involving only the derivatives in the $x_1$ and the $x_2$ directions.
It follows that
\begin{align}
	\Vert Z \nabla p\Vert_{L^{\infty}}
	 \lec
	  \Vert D^2 p\Vert_{H^2_\cco}+\Vert \nabla p\Vert_{H^2_\cco}
	  .\label{EQ15}
\end{align}
Now, combining \eqref{EQ.Pre}--\eqref{EQ15}, 
we arrive at
\begin{align}
	\Vert Zu(t)\Vert_{L^{\infty}}
	\lec
	1
	+\int_0^t \left(\Vert \nabla u\Vert_{L^{\infty}}\Vert u\Vert_{W^{1,\infty}_\cco}
	+\Vert u(t)\Vert_{W^{2,\infty}_\cco}+\Vert \nabla u(t)\Vert_{L^\infty}+1
	\right) \,ds
	,\label{EQ16}
\end{align}
for $t \in [0,T]$ and $T<T^*$.
We note that integrating \eqref{EQ09} in time and using \eqref{EQas} yields
\begin{align}
	\int_0^{T*} \Vert \nabla u\Vert_{L^{\infty}}\,ds
	\lec 1
	,\llabel{EQ17}
\end{align}
hence, \eqref{EQ16} implies 
\begin{align}
	\limsup_{T \to T^*}\sup_{[0,T]}\Vert u\Vert_{W^{1,\infty}_\cco}
	\lec
	1
	.\label{EQ18}
\end{align}
Collecting \eqref{EQ09} and \eqref{EQ18}, we also obtain
\begin{align}
	\limsup_{T \to T^*}\sup_{[0,T]}\Vert \nabla u(t)\Vert_{L^{\infty}}
	\lec
	1
	.\label{EQ19}
\end{align}

Finally, we establish that $\Vert u\Vert_{2,\infty}$ stays bounded up to time $T^*$.
Now, repeating the steps in the above paragraph for $Z^\alpha u$ and $|\alpha|=2$, we arrive at
\begin{align}
	\Vert Z^\alpha u\Vert_{L^{\infty}}
	 \lec
	  \Vert Z^\alpha u_0\Vert_{L^{\infty}}
	  +\int_0^t
	   (\Vert u\cdot \nabla Z^\alpha u - Z^\alpha(u\cdot \nabla u)\Vert_{L^\infty}
	   +\Vert Z^\alpha \nabla p\Vert_{L^\infty})\,ds
	   .\label{EQ20}
\end{align}
We may bound the commutator term by writing
\begin{align}
	\Vert u\cdot \nabla Z^\alpha u - Z^\alpha(u\cdot \nabla u)\Vert_{L^\infty}
	 \lec
	  \Vert u\Vert_{W^{2,\infty}_\cco}
	   (\Vert u\Vert_{W^{2,\infty}_\cco}
	    +\Vert u\Vert_{W^{1,\infty}_\cco}
	     +\Vert \nabla u\Vert_{L^\infty})
	     ,\llabel{EQ21}
\end{align}
while for the pressure term, we have
\begin{align}
	\Vert \nabla p\Vert_{W^{2,\infty}_\cco}
	 \lec
	  \Vert D^2 p\Vert_{H^3_\cco}+\Vert \nabla p\Vert_{H^3_\cco}
	 .\llabel{EQ22}
\end{align}
Utilizing \eqref{EQ.Pre}, \eqref{EQ10}, and \eqref{EQ18}--\eqref{EQ20},
we note that 
\begin{align}
	\limsup_{T \to T^*}\sup_{[0,T]} \Vert u\Vert_{W^{2,\infty}_\cco} \lec 1.
	\llabel{EQ23}
\end{align}

In summary, we have shown that 
\begin{align}
	\limsup_{T \to T^*} \sup_{[0,T]} 
	 (\Vert u\Vert_{H^4_\cco} + \Vert u\Vert_{W^{2,\infty}_\cco} + \Vert u\Vert_{W^{1,\infty}})
	 \le M
	 ,\llabel{EQ62}
\end{align}
where $M>0$ is an absolute constant. Now, recalling that $T>0$ given by Theorem~\ref{T01}
depends on norms of the initial datum, we may choose $\epsilon>0$ sufficiently small
and invoke Theorem~\ref{T01} to solve the Euler equations on $[T^*-\epsilon, T^*-\epsilon+T]$ with the initial datum $u(T^*-\epsilon)$, where $T-\epsilon>0$.
Then, utilizing Lemma~\ref{L03} and the uniqueness of conormal solutions, we may extend $u$ beyond $T^*$ 
and this concludes the proof of~Theorem~\ref{T02}.
\end{proof}

\section*{Acknowledgments}
\rm
The author was supported in part by the NSF grant DMS-2205493 and the USC Dornsife Summer Research Award.

\ifnum\sketches=1

\fi


\begin{thebibliography}{[GHKR]}
\small

\bibitem[A]{A} 
M.S.~Ayd{\i}n,  
\emph{Inviscid limit on $L^p$-based Sobolev conormal spaces for the 3D Navier-Stokes equations with the Navier boundary conditions},
arXiv:2502.03599.
	
\bibitem[AK1]{AK1} 
M.S.~Ayd\i n, and I.~Kukavica, 
  \emph{Euler equations in {S}obolev conormal spaces},
  Res. Math. Sci.~\textbf{12} (2025), no.~2, Paper No. 33.
    
\bibitem[AK2]{AK2} 
M.S.~Ayd{\i}n, and I.~Kukavica, 
\emph{Uniform bounds and the inviscid limit for the {N}avier-{S}tokes equations with {N}avier boundary conditions},
arXiv:2404.17111. 

\bibitem[BD]{BD} 
H.~Bahouri, and B.~Dehman, 
\emph{Remarques sur l'apparition de singularit\'es dans les
	\'ecoulements eul\'eriens incompressibles \`a{} donn\'ee
	initiale h\"old\'erienne},
J. Math. Pures Appl. (9)~\textbf{73} (1994), no.~4, 335--346.

\bibitem[BILN]{BILN} 
A.V.~Busuioc, D.~Iftimie, M.C.~Lopes~Filho, and H.J.~Nussenzveig~Lopes,
\emph{Uniform time of existence for the alpha {E}uler equations},
J. Funct. Anal.~\textbf{271} (2016), no.~5, 1341--1375.
    
\bibitem[BKM]{BKM} 
J.T.~Beale, T.~Kato, and A.~Majda, 
\emph{Remarks on the breakdown of smooth solutions for the {$3$}-{D}
	{E}uler equations},
Comm. Math. Phys.~\textbf{94} (1984), no.~1, 61--66.

\bibitem[Ch1]{Ch1} 
D.~Chae,  
\emph{On the well-posedness of the {E}uler equations in the {B}esov
	and the {T}riebel-{L}izorkin spaces},
S\=urikaisekikenky\=usho Koky\=uroku (2001), no.~1234, 42--57.

\bibitem[Ch2]{Ch2} 
D.~Chae,  
\emph{On the well-posedness of the {E}uler equations in the
	{T}riebel-{L}izorkin spaces},
Comm. Pure Appl. Math.~\textbf{55} (2002), no.~5, 654--678.

\bibitem[Ch3]{Ch3} 
D.~Chae,  
\emph{Remarks on the blow-up criterion of the three-dimensional
	{E}uler equations},
Nonlinearity~\textbf{18} (2005), no.~3, 1021--1029.

\bibitem[Ch4]{Ch4} 
D.~Chae,  
\emph{On the continuation principles for the {E}uler equations and
	the quasi-geostrophic equation},
J. Differential Equations~\textbf{2} (2006), no.~2, 640--651.

\bibitem[Ch5]{Ch5} 
D.~Chae,  
\emph{On the finite-time singularities of the 3{D} incompressible
	{E}uler equations},
Comm. Pure Appl. Math.~\textbf{60} (2007), no.~4, 597--617.

\bibitem[Co]{Co} 
P.~Constantin,  
\emph{An {E}ulerian-{L}agrangian approach for incompressible fluids:
	local theory},
J. Amer. Math. Soc.~\textbf{14} (2001), no.~2, 263--278.

\bibitem[CC1]{CC1} 
D.~Chae, and P.~Constantin,  
\emph{Remarks on type {I} blow-up for the 3{D} {E}uler equations and
	the 2{D} {B}oussinesq equations},
J. Nonlinear Sci.~\textbf{31} (2021), no.~5, Paper No. 77,16.

\bibitem[CC2]{CC2} 
D.~Chae, and P.~Constantin,  
\emph{On a type {I} singularity condition in terms of the pressure
	for the {E}uler equations in {$\Bbb{R}^3$}},
Int. Math. Res. Not. IMRN~\textbf{12} (2022), no.~12, 9013--9023.

\bibitem[CH]{CH} 
J.~Chen, and T.Y.~Hou,  
\emph{Finite time blowup of 2{D} {B}oussinesq and 3{D} {E}uler
	equations with {$C^{1,\alpha}$} velocity and boundary},
Comm. Math. Phys.~\textbf{383} (2021), no.~3, 1559--1667.

\bibitem[CFM]{CFM} 
P.~Constantin, C.~Fefferman, and A.~Majda, 
\emph{Geometric constraints on potentially singular solutions for
	the {$3$}-{D} {E}uler equations},
Comm. Partial Differential Equations~\textbf{21} (1996), no.~21, 599--571.

\bibitem[CI]{CI} 
D.~Chae, and O.Y.~Imanuvilov,  
\emph{Generic solvability of the axisymmetric {$3$}-{D} {E}uler
	equations and the {$2$}-{D} {B}oussinesq equations},
J. Differential Equations~\textbf{156} (1999), no.~1, 1--17.

\bibitem[CK]{CK} 
D.~Chae, and N.~Kim,  
\emph{On the breakdown of axisymmetric smooth solutions for the
	3-{D} {E}uler equations},
Comm. Math. Phys.~\textbf{178} (1996), no.~2, 391--398.

\bibitem[CW1]{CW1} 
D.~Chae, and J.~Wolf,  
\emph{Removing type {II} singularities off the axis for the three
	dimensional axisymmetric {E}uler equations},
Arch. Ration. Mech. Anal.~\textbf{234} (2019), no.~3, 1041--1089.

\bibitem[CW2]{CW2} 
D.~Chae, and J.~Wolf,  
\emph{Localized non blow-up criterion of the {B}eale-{K}ato-{M}ajda
	type for the 3{D} {E}uler equations},
Math. Ann.~\textbf{383} (2022), no.~3-4, 837--865.

\bibitem[CW3]{CW3} 
D.~Chae, and J.~Wolf,  
\emph{Localized blow-up criterion for {$ C^{ 1, \alpha}$} solutions
	to the 3{D} incompressible {E}uler equations},
J. Math. Fluid Mech.~\textbf{25} (2023), no.~3, Paper No. 68, 29.

\bibitem[DHY1]{DHY1} 
J.~Deng, T.Y.~Hou, and X.~Yu,  
\emph{Geometric properties and nonblowup of 3{D} incompressible
	{E}uler flow},
Comm. Partial Differential Equations~\textbf{30} (2005), no.~1-3, 225--243.

\bibitem[DHY2]{DHY2} 
J.~Deng, T.Y.~Hou, and X.~Yu,  
\emph{Improved geometric conditions for non-blowup of the 3{D}
	incompressible {E}uler equation},
Comm. Partial Differential Equations~\textbf{31} (2006), no.~1-3, 293--306.

\bibitem[E]{E} 
T.~Elgindi,  
\emph{Finite-time singularity formation for {$C^{1,\alpha}$}
	solutions to the incompressible {E}uler equations on {$\Bbb
		R^3$}},
Ann. of Math. (2)~\textbf{194} (2021), no.~3, 647--727.

\bibitem[F]{F} 
A.B.~Ferrari,   
\emph{On the blow-up of solutions of the {$3$}-{D} {E}uler equations
	in a bounded domain},
Comm. Math. Phys.~\textbf{155} (1993), no.~2, 277--294.

\bibitem[GHKR]{GHKR} 
J.D.~Gibbon, D.D.~Holm, and R.M.~Kerr, and I.~Roulstone,  
\emph{Quaternions and particle dynamics in the {E}uler fluid
	equations},
Nonlinearity~\textbf{19} (2006), no.~8, 1969--1983.

\bibitem[K]{K} 
N.~Kim,   
\emph{Criterion for blow-up in the {E}uler equations via certain
	physical quantities},
J. Korean Soc. Ind. Appl. Math.~\textbf{16} (2012), no.~4, 243--248.

\bibitem[KOT]{KOT} 
H.~Kozono, T~Ogawa, and Y.~Taniuchi,  
\emph{The critical {S}obolev inequalities in {B}esov spaces and
	regularity criterion to some semi-linear evolution equations}, 
Math. Z.~\textbf{242} (2002), no.~2, 251--278. 

\bibitem[KT]{KT} 
H.~Kozono, and Y.~Taniuchi, 
\emph{Limiting case of the {S}obolev inequality in {BMO}, with
	application to the {E}uler equations},
Comm. Math. Phys.~\textbf{214} (200), no.~1, 191--200.

\bibitem[MR1]{MR1} 
N.~Masmoudi, and F.~Rousset,  
\emph{Uniform Regularity for the Navier–Stokes Equation with Navier Boundary Condition}, 
Arch.\ Ration.\ Mech.\ Anal.~\textbf{203} (2012), no.~2, 529--575. 
	
\bibitem[MR2]{MR2} 
N.~Masmoudi, and F.~Rousset,  
\emph{Uniform regularity and vanishing viscosity limit for the free
	surface {N}avier-{S}tokes equations}, 
Arch. Ration. Mech. Anal.~\textbf{223} (2017), no.~1, 301--417. 

\bibitem[OT]{OT} 
T~Ogawa, and Y.~Taniuchi,  
\emph{On blow-up criteria of smooth solutions to the 3-{D} {E}uler
	equations in a bounded domain}, 
J. Differential Equations~\textbf{190} (2003), no.~1, 39--63. 

\bibitem[P]{P} 
G.~Ponce,  
\emph{Remarks on a paper: ``{R}emarks on the breakdown of smooth
		solutions for the {$3$}-{D} {E}uler equations'' [{C}omm.
		{M}ath. {P}hys. {\bf 94} (1984), no. 1, 61--66; {MR}0763762
		(85j:35154)] by {J}. {T}. {B}eale, {T}. {K}ato and {A}.
		{M}ajda}, 
Comm. Math. Phys.~\textbf{98} (1985), no.~3, 349--353. 

\bibitem[SY]{SY} 
T.~Shirota, and T.~Yanagisawa,
\emph{A continuation principle for the {$3$}-{D} {E}uler equations
	for incompressible fluids in a bounded domain}, 
Proc. Japan Acad. Ser. A Math. Sci.~\textbf{69} (1993), no.~3, 77--82. 

\bibitem[Z]{Z} 
W.M.~Zaj\c{a}czkowski,
\emph{Remarks on the breakdown of smooth solutions for the {$3$}-d
	{E}uler equations in a bounded domain}, 
Bull. Polish Acad. Sci. Math.~\textbf{37} (1989), no.~1-6, 169--181. 

\end{thebibliography}
\end{document}